\documentclass[10pt]{article} 

\usepackage[utf8x]{inputenc}
\usepackage{amsmath, amsfonts, mathrsfs, textcomp,amssymb,natbib,bbm,natbib}
\usepackage[breaklinks]{hyperref}
\usepackage{hypernat}
\usepackage[dvips]{graphicx}

\newtheorem{theorem}{Theorem}

\newtheorem{remark}{Remark}
\newtheorem{corollary}{Corollary}

\newtheorem{proposition}{Proposition}
\newcommand{\ninf}{\mbox{ as } N\to\infty}

\newcommand{\szer}{\mbox{ as } s\to 0}

\newcommand{\acal}{{\cal A}}

\newcommand{\bcal}{{\cal B }}
\newcommand{\ccal}{{\cal C }}
\newcommand{\dcal}{{\cal D }}

\newcommand{\norm}[1]{\left|\left|\, #1\,\right|\right|}
\newcommand{\mb}[1]{\mathbf{#1}}

\newcommand{\dt}{\delta}

\newcommand{\E}{\mathrm{E}}
\newcommand{\V}{\mathrm{Var}}

\newcommand{\beq}{ \begin{equation}}
\newcommand{\eeq}{ \end{equation}}
\newcommand{\beqr}{ \begin{eqnarray}}
\newcommand{\eeqr}{ \end{eqnarray}}
\newcommand{\beqrn}{ \begin{eqnarray*}}
\newcommand{\eeqrn}{ \end{eqnarray*}}
\newcommand{\bye}{\end{document}}

\title{Stratified Monte Carlo quadrature for continuous random fields}
\author{Konrad Abramowicz,$\quad$ Oleg Seleznjev,\\
 Department of Mathematics and Mathematical Statistics\\
Ume{\aa } University, SE-901 87 Ume{\aa}, Sweden }
        \date{\today}

\begin{document}

\maketitle

\begin{abstract}
We consider the problem of numerical approximation of integrals of random fields over a unit hypercube.
We use a stratified Monte Carlo quadrature and measure the approximation performance  by the mean squared error.
The quadrature is defined by a finite number of stratified randomly chosen observations with the partition generated by a rectangular grid (or design).
We study the class of locally stationary random fields whose local behavior is like a
fractional Brownian field in the mean square sense and find the asymptotic approximation accuracy for a sequence of designs for large number of the observations.
For the H\"{o}lder class of random functions, we provide an upper bound for the approximation error.
Additionally, for a certain class of isotropic random functions with an isolated singularity at the origin,
we construct a sequence of designs  eliminating the effect of the singularity point.
\end{abstract}


 \baselineskip=3.4 ex


\section{Introduction}

Let  $X(\mb{t}),\,\mb{t}\in[0,1]^d$, $d\geq 1$, be a continuous random field with finite second moment.
We consider the problem of numerical approximation of the integral of $X$ over the unit hypercube using finite number of observations.
The approximation accuracy is measured by the mean squared error.
We use a stratified Monte Carlo quadrature (sMCQ) for the integral approximation introduced for deterministic functions by \citet{Haber1966}.
The quadrature is defined by stratified random observations with the partition generated by a rectangular grid (or design).
We use cross regular sequences of designs, generalizing the well known regular sequences pioneered by
\citet{SacksYlvisaker1966}.
We focus on random fields satisfying a local stationarity condition proposed for stochastic processes by \citet{Berman1974} and extended for random fields in \citet{AbramowiczSeleznjev2011}.
Approximation of random functions from this class is studied in, e.g., \citet{Seleznjev2000, HuslerPiterbargSeleznjev2003,  AbramowiczSeleznjev2011, Abramowicz2011}.
For quadratic mean (q.m.) continuous locally stationary random functions, we derive an exact asymptotic behavior of the approximation accuracy.
We propose a method for the asymptotically optimal sampling point distribution between the mesh dimensions.
We also study optimality of grid allocation along coordinates and provide asymptotic optimality results in the one-dimensional case.
For q.m.\ continuous fields satisfying a H\"{o}lder type condition, we determine an upper bound for the approximation accuracy.
Furthermore, we investigate a certain class of random fields with different q.m. smoothness at the origin (isolated singularity), and construct
sequences of designs eliminating the effect of the singularity point.

Approximation of integrals of random functions is an important problem arising in many research and applied areas, like environmental and geosciences \citep{Ripley2004}, 
communication theory and signal processing \citep{Masry2009}.
Regular sampling designs for estimating integrals of stochastic processes are studied in \citet{Benhenni1992}.
%
Random designs of sampling points, including stratified sampling for stochastic processes, are investigated in \citet{Schoenfelder1982, CambanisMasry1992}.
Minimax results for estimating integrals of analytical processes are presented in \citet{BenhenniIstas1998}.
Prediction of integrals of stationary random fields using the observations on a lattice is discussed in \citet{Stein1995}.
Quadratures for smooth isotropic random functions are investigated in \citet{Ritter1997, Stein1995_2}.
Multivariate numerical integration of random fields satisfying {S}acks-{Y}lvisaker conditions is studied in \citet{Ritter1995}.
\citet{Ritter2000} contains a survey of various random function approximation and integration problems.
\citet{Novak1988} includes a detailed discussion of deterministic and Monte Carlo (randomized) linear methods in various computational problems.
We refer to \citet{Adler2007} for a comprehensive summary of the general theory of random fields.

The paper is organized as follows. First we introduce a basic notation.
In Section 2, we consider a stratified Monte Carlo quadrature for continuous random fields which local behavior is like a fractional Brownian field in the mean square sense.
We derive an exact asymptotics and a formula for the optimal interdimensional sampling point distribution.
Further, we provide an upper bound for the approximation accuracy for q.m.\ continuous fields satisfying H\"older type conditions.
In the second part of this section, we study random fields with an isolated singularity at the origin and construct sequences of designs eliminating the effect of the singularity.
In Section 3, we present the results of numerical experiments, while Section 4 contains the proofs of the statements from Section 2.

\subsection{Basic notation}

Let $X=X(\mathbf{t}),\textbf{t}\in\dcal:= [0,1]^d$, $d\geq1$, be a random field defined on a probability space $(\Omega,\mathscr{F},P)$. Assume that for every $\textbf{t}$, the random
variable $X(\textbf{t})$ lies in the normed linear space $L^2(\Omega)=L^{2}(\Omega,\mathscr{F},P)$ of
 random variables with finite second moment and identified equivalent elements with respect to $P$.
We set $||\xi||:=\left(\mathrm{E}\xi^2\right)^{1/2}$ for all $\xi \in L^2(\Omega)$.
We are interested in a numerical approximation of
$$
I(X)=\int_\dcal X(\mb{t})d\mb{t}
$$
by a quadrature based on $N$ observations for random fields from a space $\ccal(\dcal)$ of q.m. continuous random fields.

We introduce the classes of random fields used throughout this paper.
For $k\le d$, let $\mathbf{l}=(l_1,\ldots,l_k)$ be a vector of positive integers such that $\sum_{j=1}^{k}l_j=d$, and let $L_i:=\sum_{j=1}^{i}l_j, i=0,\ldots,k$, $L_0=0$, be the
sequence of its cumulative sums. Then the vector $\mathbf{l}$ defines the \mbox{\textit{l-decomposition}} of
$\dcal$ into $\dcal^1\times\ldots\times\dcal^k$, with the $l_j$-cube $\dcal^j=[0,1]^{l_j}$, $j=1,\ldots,k$.
For any $\mathbf{s}\in\dcal$, we denote by $\mathbf{s}^j$ the coordinates vector corresponding to the $j$-th component
of the decomposition, i.e.,
$$
 \mathbf{s}^j=\mathbf{s}^j(\mb{l}):=(s_{L_{j-1}+1},\ldots, s_{L_{j}}) \in \dcal^j,\quad j=1,\ldots,k.
$$
For a vector $\boldsymbol\alpha=(\alpha_1,\ldots,\alpha_k)$, $0<\alpha_j<2$, $j=1,\ldots,k$, and the decomposition vector $\mathbf{l}=(l_1,\ldots,l_k)$,
let
$$
\norm{\mathbf{s}}_{\boldsymbol\alpha}:=\sum_{j=1}^{k}\norm{\mathbf{s}^j}^{\alpha_j} \quad \mbox{ for all } \mathbf{s}\in\dcal
$$
with the Euclidean norms $||\mathbf{s}^j||, j=1,\ldots,k$.\\

\noindent  For a hyperrectangle $\acal=[a_1,b_1]\times\ldots\times[a_d,b_d]\subset\dcal$ and a random field $X\in\ccal(\acal)$, we say that\\
(i) $X\in\ccal_\mathbf{l}^{\boldsymbol\alpha}(\acal,C)$ if for some $\boldsymbol\alpha$,  $\mathbf{l}$,  and a positive constant $C$,
the random field $X$ satisfies the H\"older condition,  i.e.,
\beq \label{def:hcont}
    \norm{X(\mathbf{t+s})-X(\mathbf{t})}^2 \leq C \norm{\mathbf{s}}_{\boldsymbol\alpha} \qquad\mbox{ for all } \mathbf{t},\mathbf{t+s}\in\acal;
\eeq
(ii) $X\in\bcal_\mathbf{l}^{\boldsymbol\alpha}(\acal,c(\cdot))$ if for some $\boldsymbol\alpha$, $\mathbf{l}$, and a vector function $c(\mb{t})=(c_1(\mb{t}),\ldots,c_k(\mb{t}))$, $\mb{t}\in\acal$, the random field $X$ is \textit{locally stationary}, i.e.,
\beq \label{def:locstat}
  \frac{\norm{X(\mathbf{t+s})-X(\mathbf{t})}^2}{\sum_{j=1}^k c_k(\mathbf{t})\norm{\mathbf{s}^j}^{\alpha_j}}\rightarrow 1\quad\mbox{as }\mathbf{s}\rightarrow 0 \mbox{ uniformly in }\mathbf{t}\in\acal,
\eeq
with positive and continuous functions $c_1(\cdot),\ldots,c_k(\cdot)$.
We assume additionally that for $j=1,\ldots,k$, the function $c_j(\cdot)$ is invariant with respect to coordinates permutation within the $j$-th component. \\\medskip

\noindent  For the classes $\ccal_\mathbf{l}^{\boldsymbol\alpha}$ and $\bcal_\mathbf{l}^{\boldsymbol\alpha}$, the withincomponent smoothness is defined by the vector $\boldsymbol\alpha=(\alpha_1,\ldots,\alpha_k)$. We denote the vector describing the smoothness for each coordinate by $\boldsymbol\alpha^*=(\alpha_1^*,\ldots,\alpha_d^*)$, where $\alpha_i^* =\alpha_j$, $i=L_{j-1}+1,\ldots,L_j$, $j=1,\ldots,k$.
Moreover, for one component fields, i.e., $k=1$ and $\boldsymbol{\alpha}=\alpha$, the corresponding H\"{o}lder and local stationary classes are denoted by $\ccal^\alpha_d$ and $\bcal^\alpha_d$, respectively.\medskip

\noindent\textbf{Example 1.}
Let $\mb{m}=(m_1,\ldots,m_k)$ be a decomposition vector of $[0,1]^m$, and $m=\sum_{j=1}^k m_j$.
Denote by $B_{\boldsymbol\beta,\mb{m}}(\mathbf{t})$, $\mathbf{t}\in[0,1]^m$,
$\boldsymbol\beta=(\beta_1,\ldots,\beta_k)$, $0<\beta_j<2$, $j=1,\ldots,k$,
an $m$-dimensional fractional Brownian field with covariance function
$r(\mathbf{t},\mathbf{s})=\frac{1}{2}\left(||\mathbf{t}||_{\boldsymbol\beta}+||\mathbf{s}||_{\boldsymbol\beta}-||\mathbf{t}-\mathbf{s}||_{\boldsymbol\beta}\right)$.
Then $B_{\boldsymbol\beta,\mb{m}}$ has stationary increments,
$$
||B_{\boldsymbol\beta,\mathbf{m}}(\mathbf{t+s})-B_{\boldsymbol\beta,\mathbf{m}}(\mathbf{t})||^2=||\mathbf{s}||_{\boldsymbol\beta}, \quad \mathbf{t}, \mathbf{t+s} \in[0,1]^m,
$$
and therefore, $B_{\boldsymbol\beta,\mb{m}}\in\bcal_\mathbf{m}^{\boldsymbol\beta}(\dcal,c(\cdot))$ with local stationarity functions $c_1(\mathbf{t})=\ldots=c_k(\mb{t})=1$, $\mathbf{t}\in[0,1]^m$.
In particular, if $k=1$, then $B_{\beta,m}(\mathbf{t}),\,\mb{t}\in[0,1]^m$, $0<\beta<2$, $m\in\mathbbm{N}$, is an $m$-dimensional fractional Brownian field with covariance function
\begin{equation}\label{eq:FBFdef}
r(\mathbf{t},\mathbf{s})=\frac{1}{2}\left(||\mathbf{t}||^{\beta}+||\mathbf{s}||^{\beta}-||\mathbf{t}-\mathbf{s}||^{\beta}\right),  \quad  \mathbf{t}, \mathbf{t+s} \in[0,1]^m.
\end{equation}\medskip\\
Let the hypercube $\dcal$ be partitioned into hyperrectangular strata by design points $T_N$, for $N\geq 1$.
We consider \textit{cross regular sequences} of grid designs \citep[see, e.g.,][]{AbramowiczSeleznjev2011}. The designs $T_N:=\{\mathbf{t}_{\mathbf{i}}=(t_{1,i_1},\ldots,t_{d,i_d}): \mathbf{i}=(i_1,\ldots,i_d)$, $0\leq i_k\leq n_k^*$, $k=1,\ldots,d\}$ are defined by the one-dimensional grids
$$
\int_0^{t_{j,i}}h_j^*(v)dv=\frac{i}{n_{j}^{*}},\quad i=0,1,\ldots,n_j^*,\quad j=1,\ldots,d,
$$
where $h^*_j(s)$, $s\in[0,1]$, $j=1,\ldots,d$,  are positive and continuous density functions, say, \textit{withindimensional densities}, and
let
$$
{h}^*(\mb{t}):=(h_1^*(t_1),\ldots,h_d^*(t_d)).
$$
The \textit{interdimensional grid distribution} of sampling points is determined by a vector function
$\pi:\mathbbm N\to\mathbbm{N}^d$:
$$
\pi^*(N):=(n_1^*(N),\ldots,n_d^*(N)),
$$
where $\lim_{N\to \infty} n_j^*(N)=\infty$, $j=1,\ldots,d$, and the condition
$$\prod_{j=1}^{d}n_j^*(N)=N$$
is satisfied. We suppress the argument $N$ for $n_j^*=n_j^*(N)$, $j=1,\ldots,d$, when doing so causes no confusion.

The introduced classes of random fields have the same smoothness and local behavior for each coordinate of the components generated by a decomposition vector $\mathbf{l}$.
Therefore we use designs with the same within- and interdimensional grid distributions within the components.
Formally, for the partition generated by a vector $\mb{l}=(l_1,\ldots,l_k)$, we consider cross regular designs $T_N$, defined by functions
$h=(h_1,\ldots,h_k)$ and $\pi(N)=(n_1(N),\ldots,n_k(N))$, in the following way:
$$
h_i^*(\cdot)\equiv h_j(\cdot),\quad
n_i^* = n_j,\quad i=L_{j-1}+1,\ldots,L_j,\quad j=1,\ldots,k.
$$
We call functions $h_1(\cdot),\ldots,h_k(\cdot)$ and $\pi(N)$  \textit{withincomponent densities} and  \textit{intercomponent grid distribution}, respectively.
The corresponding property of a design $T_N$ is denoted by: $T_N$ is $cRS(h,\pi,\mb{l})$.
If $d=1$, then $\mb{l}=1$, $\pi(N)=\pi_1(N)=N$, and
the cross regular sequences become regular sequences introduced by \citet{SacksYlvisaker1966}. We denote such property of the design by: $T_N$ is $RS(h)$.\medskip\\
For a given cross regular grid  design, the hypercube $\dcal$ is partitioned
into $N$ disjoint {hyperrectangular strata}  $\dcal_{\mathbf{i}}$, $\mb{i}\in\mb{I}$, where
$\mb{I}:=\{\mathbf{i}=(i_1,\ldots,i_d)$, $0\leq i_k\leq n_k^*-1$, $k=1,\ldots,d\}$. Let $\mb{1}_d=(1,\ldots,1)$ and $\mb{0}_d=(0,\ldots,0)$ denote a $d$-dimensional vectors of ones and zeros, respectively.
The hyperrectangle $\dcal_{\mathbf{i}}$ is determined by the vertex
\mbox{$\mathbf{t_{i}}=(t_{1,i_1},\ldots,t_{d,i_d})$}
and the main diagonal $\mathbf{r}_{\mathbf{i}}:=\mathbf{t}_{\mb{i}+\mb{1}_d}-\mathbf{t_{i}}$, i.e.,
$$
\dcal_{\mathbf{i}}:=\left\{\mathbf{t}: \mathbf{t}=\mathbf{t}_{\mathbf{i}}+\mathbf{r}_{\mathbf{i}} *\mathbf{s}, \mathbf{s}\in[0,1]^d\right\},
$$
where $'*'$ denotes the coordinatewise multiplication, i.e., for $\mathbf{x}=(x_1,\ldots,x_d)$ and $\mathbf{y}=(y_1,\ldots,y_d)$,
{$\mathbf{x}\ast\mathbf{y}:=(x_1y_1,\ldots,x_dy_d)$}.\medskip\\
Let $|\dcal_\mb{i}|$ denote the volume of the hyperrectangle $\dcal_\mb{i}$. For a random field $X\in{\ccal}(\dcal)$, we define a \textit{stratified Monte Carlo quadrature}
(sMCQ) on a partition generated by $T_N$
$$
I_N(X,T_N):=I_N(X,T_N(h,\pi,\mb{l})) = \sum_{\mb{i}\in \mb{I}} X(\boldsymbol\eta_\mb{i})|\dcal_\mb{i}|,
$$
where $\boldsymbol\eta_\mb{i}$ is uniformly distributed in the stratum $\dcal_\mb{i}$, $\mb{i}\in\mb{I}$.
Such defined quadrature is a modification of a well known midpoint quadrature.

\section{Results}
Let $B_{\beta,m}(\mathbf{t}),\mathbf{t}\in\mathbb{R}^{m}_{+}$, $0<\beta<2$, $m\in\mathbb{N}$, denote an $m$-dimensional fractional
Brownian field with covariance function \eqref{eq:FBFdef}. For any $\mb{u}\in\mathbb{R}^{m}_{+}$, we denote
\begin{equation}\label{def:bconst}
b_{\beta,m}(\mb{u})
           =\frac{1}{2}\int_{[0,1]^m}\int_{[0,1]^m}\norm{\mb{u}\ast(\mb{t}-\mb{v})}^{\beta}d\mb{t}d\mb{v}
           =\E \left( \int_{[0,1]^m} B_{\beta,m}(\mb{u}\ast\mb{t})d\mb{t} - B_{\beta,m}(\mb{u}\ast{\boldsymbol\eta}) \right)^2,
\end{equation}
where $\boldsymbol\eta$ is uniformly distributed in the unit $m$-hypercube. Then $b_{\beta,m}(\mb{u})$ corresponds to
the mean squared error (MSE) of a sMCQ based on one observation for a field $B_{\beta,m}(\mb{u}\ast\mb{t})$, $\mb{t}\in[0,1]^m$.

In the following theorem, we provide an exact asymptotics for the accuracy of a sMCQ for
locally stationary random fields when cross regular sequences of grid designs are used.
\begin{theorem}\label{Th:Main}
Let $X\in\bcal^{\boldsymbol\alpha}_{\mb{l}}(\dcal,c(\cdot))$ be a random field and let $I(X)$ be approximated by sMCQ $I_N(X,T_N)$, where $T_N$ is $cRS(h,\pi,\mb{l})$. Then
$$
\norm{I(X)-I_N(X,T_N)}^2\sim \frac{1}{N} \sum_{j=1}^k \frac{v_j}{n_j^{\alpha_j}}\mbox{ as }N \to \infty,
$$
where
$$
v_j:=\int_\dcal c_j(\mb{t}) b_{\alpha_j,l_j}(D_j(\mb{t}^j)) \prod_{m=1}^{d} {h^*_m(t_m)}^{-1}d\mb{t}>0
$$
and $D_j(\mb{t}^j):=(1/h_j(t_{L_{j-1}+1}),\ldots,1/h_j(t_{L_j}))$.
\end{theorem}

\begin{remark}
If $T_N$ is a systematic sampling, i.e., all withincomponent grid distributions are uniform, \mbox{$h_j(s)=1$}, $s\in[0,1]$, $j=1,\ldots,k$, then
the asymptotic constants are reduced to
$$
v_j=\tilde b_{\alpha_j,l_j} \int_{\dcal} c_j(\mathbf{t})d\mathbf{t},\qquad j=1,\ldots,k,
$$
where $\tilde b_{\alpha_j,l_j}:=b_{\alpha_j,l_j}(\mathbf{1}_{l_j})$.
\end{remark}

The next theorem presents an asymptotically optimal intercomponent grid distribution for a given total number of sampling points $N$.
We define
$$
\rho:=\left(\sum_{i=1}^{k}\frac{l_i}{\alpha_i}\right)^{-1}\!\!\!\!=\left(\sum_{i=1}^{d}\frac{1}{\alpha_i^*}\right)^{-1}\!\!\!\!, \qquad \kappa:=\prod_{j=1}^{k}v_{j}^{l_j/\alpha_j},
$$
where $d\!\cdot\!\rho$ is the harmonic mean of the smoothness parameters $\alpha_j^*,\,j=1,\ldots,d$.
\begin{theorem}\label{Th:DimOpt}
Let $X\in\bcal^{\boldsymbol\alpha}_{\mb{l}}(\dcal,c(\cdot))$ be a random field and let $I(X)$ be approximated by sMCQ $I_N(X,T_N)$, where $T_N$ is $cRS(h,\pi,\mb{l})$. Then
\beq\label{eq:Th2}
\norm{I(X)-I_N(X,T_N)}^2 \gtrsim k\ \frac{\kappa^{\rho}}{N^{1+\rho}} \mbox{ as } N\to \infty.
\eeq
Moreover, for the asymptotically optimal intercomponent grid allocation,
\beq\label{opt_pi}
n_{j,opt}\sim \frac{\ {v_{j}^{1/\alpha_j}}}{\kappa^{\rho/\alpha_j}}N^{\rho/\alpha_j} \mbox{ as } N\to \infty, \quad j=1,\ldots,k,
\eeq
the equality in \eqref{eq:Th2} is attained asymptotically.
\end{theorem}

In a general setting, numerical
procedures can be used for finding optimal densities. However, in practice such methods are very computationally demanding.
We present a simplification of the asymptotic constant expression for one-dimensional components. For a random field $X\in\bcal_\mb{l}^{\boldsymbol\alpha}(\dcal,c(\cdot))$ define
$$
\begin{aligned}
Q_j(t_{L_j})&:=\int_{[0,1]^{d-1}}c_j(\mathbf{t}) \prod_{\substack{m=1\\ m\neq L_j}}^{d}{h^*_m(t_m)}^{-1} dt_1 \ldots dt_{{L_j-1}}dt_{{L_j}+1}\ldots dt_{d},\quad &t_{L_j}\in[0,1],\quad j=1,\ldots,k.\\
\end{aligned}
$$
Moreover, for $0<\beta<2$, let
$$
a_{\beta}:=\frac{1}{(1+\beta)(2+\beta)}.
$$

\begin{proposition}\label{Prop:OneDim} Let $X\in\bcal^{\boldsymbol\alpha}_{\mb{l}}(\dcal,c(\cdot))$ be a random field and let $I(X)$ be approximated by sMCQ $I_N(X,T_N)$, where $T_N$ is $cRS(h,\pi,\mb{l})$.
 If for some $j$, $1\leq j\leq k$,  $l_j=1$, then for any regular density $h_j(\cdot)$, we have
 \begin{equation*}
 v_j=a_{\alpha_j}\int_0^1 Q_j(t) h_j(t)^{-(1+\alpha_j)} dt.
 \end{equation*}
The $j$-th withincomponent density minimizing $v_j$ is given by
\begin{equation*}
h_{j,opt}(t)=\frac{Q_j(t)^{\gamma_j}}{\int_0^1 Q_j(\tau)^{\gamma_j}d\tau}, \qquad t\in[0,1],
\end{equation*}
where $\gamma_j:=1/(2+\alpha_j)$. Furthermore, for such density, we get $$v_{j,opt}=a_{\alpha_j}\left(\int_0^1 Q_j(t)^{\gamma_j}dt\right)^{1/\gamma_j}.$$
\end{proposition}

As a direct implication of Proposition \ref{Prop:OneDim}, we obtain the following asymptotic result for the approximation of integral of locally stationary
stochastic processes by a sMCQ, with regular sequences of grid designs. Further, in this case, we get the exact formula for the density minimizing the asymptotic constant.

\begin{corollary}\label{Cor:MainOneDim}
Let $X\in\bcal^{\alpha}_1([0,1],c(\cdot))$ be a random process and let $I(X)$ be approximated by sMCQ $I_N(X,T_N)$, where $T_N$ is $RS(h)$. Then
$$
\lim_{N\to\infty} N^{1+\alpha}\norm{I(X)-I_N(X,T_N)}^2 = a_\alpha \int_0^1 c(t)h(t)^{-(1+\alpha)}dt.
$$
The density minimizing the asymptotic constant is given by
\begin{equation}\label{eq:OneDimOpt}
h_{opt}(t)=\frac{c(t)^{\gamma}}{\int_0^1 c(\tau)^{\gamma}d\tau}, \qquad t\in[0,1],
\end{equation}
where $\gamma:=1/(2+\alpha)$. Furthermore, for such density, we get
$$
\lim_{N\to\infty} N^{1+\alpha}\norm{I(X)-I_N(X,T_N)}^2 = a_{\alpha}\left(\int_0^1 c(t)^{\gamma}dt\right)^{1/\gamma}.
$$
\end{corollary}\bigskip

Now we focus on random fields satisfying the introduced H\"{o}lder type condition.
The following proposition provides an upper bound for the accuracy of sMCQ for H\"{o}lder classes of continuous fields.
In addition, we present the intercomponent grid distribution leading to an increased rate of the upper bound.

\begin{proposition}\label{Prop:Holder}
Let $X\in\ccal^{\boldsymbol\alpha}_{\mb{l}}(\dcal,C)$ be a random field and let $I(X)$ be approximated by sMCQ $I_N(X,T_N)$, where $T_N$ is $cRS(h,\pi,\mb{l})$. Then
\begin{equation}\label{eq:HoldType1}
\norm{I(X)-I_N(X,T_N)}^2\leq \frac{C}{N}\sum_{j=1}^{k}\frac{d_j}{n^{\alpha_j}}
\end{equation}
for positive constants $d_1,\ldots,d_k$. Moreover if $n_j\sim N^{\rho/\alpha_j}$, $j=1,\ldots,k$, then
$$
\norm{I(X)-I_N(X,T_N)}^2=\mathrm{O}\left(N^{-(1+\rho)}\right)\ninf.
$$
\end{proposition}
The approximation rates obtained in the above proposition are optimal in a certain sense, i.e., the rate of
convergence can not be improved in general for random functions satisfying H\"{o}lder type condition \citep[see, e.g., ][]{Ritter2000}.
The rate of the upper bound corresponds to the optimal rate of Monte Carlo methods for the anisotropic H\"{o}lder-Nikolskii class, which is
a deterministic analogue of the introduced H\"{o}lder class \citep[see, e.g.,][]{Peixin2005}.

\begin{remark} It follows from the proof of Proposition \ref{Prop:Holder} that \eqref{eq:HoldType1} holds if
$$
d_j=a_{\alpha_j} l_j^{1+\alpha_j/2} C_j^{\alpha_j}\prod_{i=1}^k C_i^{l_i},\quad j=1,\ldots,k,
$$
where $C_j:=1/\min_{s\in[0,1]}h_j(s)$, $j=1,\ldots,k$. Therefore the constants depend only on the parameters of the H\"{o}lder class and the corresponding sampling design.
\end{remark}

\subsection{Point singularity at the origin}
In this subsection, we focus on one component random fields, i.e., $k=1$, $\mb{l}=d$, $\boldsymbol\alpha=\alpha$,
and consider the case of an isolated point singularity at the origin.
More precisely, let a random function $X(\mb{t})$, $t\in[0,1]^d$, satisfy the smoothness condition \eqref{def:hcont} with $\boldsymbol\alpha=\beta$, $\beta\in(0,2)$, for $\mb{t}\in[0,1]^d$.
In addition, let $X$ be locally stationary, \eqref{def:locstat}, with parameter $\alpha>\beta$,  on any hyperrectangle $\acal \subset [0,1]^d\backslash\{\mb{0}_d\}$.
We construct sequences of grid designs with an asymptotic approximation rate $N^{-(1+\alpha/d)}$.\medskip\\
The definition of $cRS$ for $k=1$ gives that $n_j=N^{1/d}$ and $h_j^*(\cdot)=h(\cdot)$, $j=1,\ldots,d$, for a positive and continuous density $h(t)$, $t\in[0,1]$.
For the density $h(\cdot)$, we define the related distribution functions
$$
H(t):=\int_0^{t}h(u)du,\qquad G(t):=H^{-1}(t)=\int_0^{t}g(v)dv,\quad t\in[0,1],
$$
i.e., $G(\cdot)$ is a quantile function for the distribution $H$. Moreover, by
\begin{equation}\label{def:qdens}
g(t):=G'(t)={1}/{h(G(t))},\quad t\in[0,1],
\end{equation}
we denote the \emph{quantile density function}.

\noindent To formulate the forthcoming results, we introduce  additional classes of random functions. For a random function $X\in\ccal(\dcal)$, we say that:\\
(iii) $X\in\ccal^{\alpha}_d(\acal, V(\cdot))$ for a hyperrectangle $\acal\subset\dcal$ if $X$ is \textit{locally H\"{o}lder continuous}, i.e.,  if
for all $\mb{t},\mb{t+s}\in\acal$,
\beq \label{def:locHol}
    \norm{X(\mb{t+s})-X(\mb{t})}^2 \leq V(\bar{\mb{t}}) \norm{\mb{s}}^\alpha,   0<\alpha< 2,
\eeq
for a positive continuous function $V(\mb{t}), \mb{t}\in \acal$, and some
$\bar{\mb{t}}\in \{\bar{\mb{t}}:\bar{\mb{t}}=\mb{t}+\mb{s}\ast\mb{u},\mb{u}\in[0,1]^d\}$.
In particular, if $V(\mb{t})=C,\,\mb{t}\in\acal$, where $C$ is a positive constant, then $X$ is {H\"{o}lder continuous};
\\
(iv)
$X\in \ccal\bcal^{\alpha}_d([0,1]^d\backslash\{\mb{0}_d\},c(\cdot),V(\cdot))$
if there exist $0<\alpha< 2$ and positive continuous functions $c(\mb{t}),V(\mb{t})$, $\mb{t}\in [0,1]^d\backslash\{\mb{0}_d\}$, such that
$X\in \ccal^{\alpha}_d(\acal,V(\cdot)) \cap \bcal^{\alpha}_d(\acal,c(\cdot))$ for any hyperrectangle $\acal~\subset~[0,1]^d\backslash\{\mb{0}_d\}$.
By definition, we have that $V(\mb{t})\ge c(\mb{t})$, $\mb{t}\in [0,1]^d\backslash\{\mb{0}_d\}$.\medskip\\
\noindent\textbf{Example 2.} Consider a zero mean random field $X_\alpha (\mb{t})$, $0<\alpha<2$, $\mb{t}\in[0,1]^d$, $d\geq 1$,  with covariance function
$r(\mb{t},\mb{s})=\exp\left(-\norm{\mb{t}-\mb{s}}^\alpha\right)$.
Let $Y_{\alpha,\beta}(\mb{t})=\norm{\mb{t}}^{\beta/2} X_\alpha(\mb{t})$, $\mb{t}\in[0,1]^d$, where $0<\beta<\alpha$.
Then
$$
\norm{Y_{\alpha,\beta}(\mb{t}+\mb{s})-Y_{\alpha,\beta}(\mb{t})}^2=\left(\norm{\mb{t}+\mb{s}}^{\beta/2}- \norm{\mb{t}}^{\beta/2} \right)^2
    + 2 \norm{\mb{t}}^{\beta/2} \norm{\mb{t}+\mb{s}}^{\beta/2} \left(1-e^{-\norm{\mb{s}}^\alpha}\right)
$$
and it follows by calculus that
 $ Y_{\alpha,\beta}\in\ccal^{\beta}_d([0,1]^d,M)\cap\ccal\bcal^{\alpha}_{d}([0,1]^d\backslash\{\mb{0}_d\},c(\cdot),V(\cdot))$
with $M=3$, $c(\mb{t})=2\norm{\mb{t}}^{\beta}$, and $V(\mb{t})=\beta^2/4\norm{\mb{t}}^{\beta-2}+2$.\medskip\\

\noindent We say that a positive function $f(\mb{t})$, $\mb{t}\in\mathbbm{R}^d$, satisfies a \textit{shifting condition} if
there exist positive constants $C_L<C_U$, $C$, and $a$ such that
\begin{equation}\label{eq:shift}
f(\mb{s})\leq C f(\mb{v})\quad\mbox{for all }\mb{s},\mb{v}\mbox{ such that } C_L\leq\frac{\norm{\mb{s}}}{\norm{\mb{v}}}\leq C_U,\quad \mb{s},\mb{v}\in[0,a]^d\backslash
\{\mathbf{0}_d\}.
\end{equation}
An example of such function is $f(\mb{t})=\norm{\mb{t}}^\alpha$ for any $0<C_L<C_U<\infty$ and $\alpha\in\mathbbm{R}$. In the one-dimensional case, the condition \eqref{eq:shift} is satisfied, e.g.,
for any function $f(\cdot)$ which is regularly varying (on the right) at the origin \citep[cf.][]{Abramowicz2011}.\bigskip

\noindent Let $X\in\ccal^{\beta}_d([0,1]^d,M)\cap\ccal\bcal^{\alpha}_{d}([0,1]^d\backslash\{\mb{0}_d\},c(\cdot),V(\cdot))$, $0<\beta<\alpha<2$.
For $\beta>\alpha-d$, we prove that under some condition on a local H\"{o}lder function $V(\cdot)$, the cross regular sequences attain the optimal approximation rate $N^{-(1+\alpha/d)}$.
Observe that $\beta>\alpha-d$ holds for all $\alpha,\beta\in(0,2)$ if $d\geq 2$ and for $d=1$ if $\beta>\alpha-1$.
Define $\mb{H}(\mb{t}):=(H(t_1),\ldots,H(t_d))$, $\mb{t}\in[0,1]^d$, and $\mb{G}(\mb{t})=:(G(t_1),\ldots,G(t_d))$, $\mb{t}\in[0,1]^d$. We formulate the following condition:\medskip\\
(C) Let  $V(\mb{G}(\cdot))$ be bounded from above by a function $R(\cdot)$ satisfying the shifting condition \eqref{eq:shift} with $C_L=1/\sqrt{3+d}$, $C_U=\sqrt{3+d}$, and such that
$\int_{[0,1]^d}R(\mb{H}(\mb{t}))d\mb{t}<\infty$.

\begin{theorem}\label{Th:RecoveryField}
Let $X\in\ccal^{\beta}_d([0,1]^d,M)\cap\ccal\bcal^{\alpha}_{d}([0,1]^d\backslash\{\mb{0}_d\},c(\cdot),V(\cdot))$, $\alpha-d<\beta<\alpha$, be a random field and let $I(X)$ be approxi\-mated by sMCQ $I_N(X,T_N)$, where $T_N$ is $cRS(h,\pi,d)$. If the local H\"{o}lder function $V(\cdot)$ satisfies the condition \textnormal{(C)}, then
\begin{equation}\label{eq:ThRecoveryField}
 \norm{I(X)-I_N(X,T_N)}^2\sim \frac{1}{N^{1+\alpha/d}}  \int_\dcal c(\mb{t}) b_{\alpha,d}(D_1(\mb{t})) \prod_{m=1}^{d} {h(t_m)}^{-1}d\mb{t} \mbox{ as }N \to \infty,
\end{equation}
where $D_1(\mb{t})=(1/h(t_1),\ldots,1/h(t_{d}))$.
\end{theorem}

\noindent Now we consider the case $d=1$ and $0<\beta\leq\alpha-1$, which is not included in the above theorem.
We consider {\it quasi regular sequences} (qRS) of sampling designs $T_N=T_N(h)$ \citep[see, e.g.,][]{Abramowicz2011}, which are a simple modification of the regular sequences.
We assume that $h(t)$ is continuous for $t\in(0,1]$, and allow it to be unbounded in $t=0$. If $h(t)$ is  unbounded in $t=0$, then $ h(t)\to +\infty$ as $t\to 0+$.
We denote this property of $T_N$ by:  $T_N$ is qRS$(h)$.
The corresponding quantile density function $g(t)$ is assumed to be continuous for $t\in[0,1]$ with the convention that $g(0)=0$ if $h(t)\to +\infty$ as $t\to 0+$.

\noindent Let $X\in\ccal^{\beta}_1([0,1],M)\cap\ccal\bcal^{\alpha}_1((0,1],c(\cdot),V(\cdot))$, $0<\beta\leq\alpha-1$. We modify the condition (C) and formulate the following condition for a local H\"{o}lder function $V(\cdot)$ and a grid generating density $h(\cdot)$:\medskip\\
(C$'$) Let  $V(G(\cdot))$ and $g(\cdot)$ be bounded from above by functions $R(\cdot)$ and $r(\cdot)$, respectively, such that
$R(\cdot)$ and $r(\cdot)$ satisfy the shifting condition \eqref{eq:shift} with $C_L=1/2$, $C_U=2$. Moreover, let $\int_0^1 R(H(t))r(H(t))^{1+\alpha}dt<\infty$, and
\beq \label{th:RvMain}
    G(s)=\mathrm{o}\left(s^{(1+\alpha)/(2+\beta)}\right) \szer.
\eeq\medskip\\
In the following theorem, we describe the class of generating densities eliminating the effect of the singularity point for the asymptotic integral approximation accuracy.
\begin{theorem}\label{Th:Recovery}
Let $X\in\ccal^{\beta}_1([0,1],M)\cap\ccal\bcal^{\alpha}_1((0,1],c(\cdot),V(\cdot))$, $0<\beta\leq\alpha-1$, be a random process and let $I(X)$ be approximated by sMCQ $I_N(X,T_N)$, where $T_N$ is $qRS(h)$.
Let for the density $h(\cdot)$ and local H\"{o}lder function $V(\cdot)$, the condition \textnormal{(C$'$)} hold. Then
\begin{equation}\label{eq:ThRecovery}
\lim_{N\to\infty} N^{1+\alpha}\norm{I(X)-I_N(X,T_N)}^2 = a_\alpha \int_0^1 c(t)h(t)^{-(1+\alpha)}dt.
\end{equation}
\end{theorem}

\begin{remark}
For $d=1$, as indicated in Corollary \ref{Cor:MainOneDim}, the density $h_{opt}(\cdot)$ minimizing the asymptotic constant in \eqref{eq:ThRecoveryField} and \eqref{eq:ThRecovery} is given by
\eqref{eq:OneDimOpt}. Thus if the condition \textnormal{(C$'$)} holds for $X$ and $h_{opt}(\cdot)$, then $h_{opt}(\cdot)$ is the asymptotically optimal density.
\end{remark}

\section{Numerical Experiments}
In this section, we present some examples illustrating the obtained results. For given withindimensional densities, \mbox{interdimensional} distributions,
and covariance functions, we use numerical integration to evaluate the mean squared error. Denote by
$$
e_N^2(X,h,\pi,\mb{l}):=\E(I(X)-I_N(X,T_N)(X,T_N(h,\pi,\mb{l})))^2
$$
the mean squared error of sMCQ $I_N(X)$ with strata generated by the grid $T_N$.
We write $h_{uni}(\cdot)$ to denote the vector of withincomponent uniform densities. Analogously, by $\pi_{uni}(\cdot)$ we denote the uniform interdimensional grid distribution, i.e.,
$n_1=\ldots=n_k$.

\begin{figure}[htb]
\begin{center}
\includegraphics[height=1.8in]{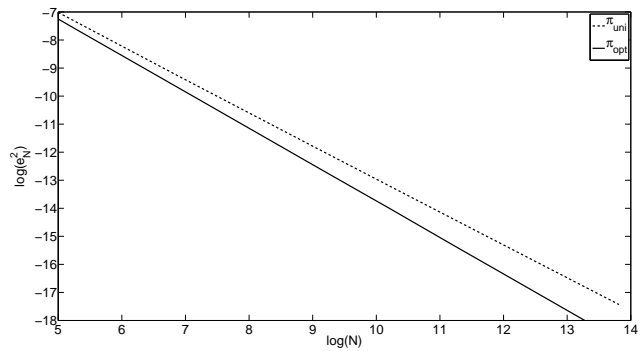}
\end{center}
\caption{The (fitted) plots of $e_N^2(\pi_{uni})$ (dashed line) and $e_N^2(\pi_{opt})$ (solid line) versus $N$ in a log-log scale.}
\label{Fig:Ex1}
\end{figure}

\noindent\textbf{Example 3.}
Let $\dcal=[0,1]^3$ and
$$
X(\mb{t})=B_{\boldsymbol\alpha,\mb{l}}(\mathbf{t}),\quad \mb{t}\in[0,1]^3,
$$
where $\boldsymbol\alpha=(3/2,1/2)$ and $\mathbf{l}=(2,1)$. Then
$X\in\bcal_{\mb{l}}^{\boldsymbol\alpha}([0,1]^3,c(\cdot))$, with $c(\mb{t})=(1,1),\,\mb{t}\in[0,1]^3$, $k=2$, $\boldsymbol\alpha^*=(3/2,3/2,1/2)$.
We compare behavior of $e_N(\pi_{uni})=e_N^2(X,h_{uni},\pi_{uni},\mb{l})$ and $e_N(\pi_{opt})=e_N^2(X,h_{uni},\pi_{opt},\mb{l})$, where the asymptotically optimal grid distribution $\pi_{opt}$ is given by Theorem~\ref{Th:DimOpt}.
Figure \ref{Fig:Ex1} shows the (fitted) plots of the mean squared errors $e_N^2(\pi_{uni})$ (dashed line) and $e_N^2(\pi_{opt})$ versus $N$ (in a log-log scale).
These plots
correspond to the following asymptotic behavior:
$$
\begin{array}{rcll}
 e_N^2(\pi_{uni}) &\sim&    C_1\, N^{-7/6}+C_2\,N^{-3/2}\,\sim\,C_1\,N^{-7/6},&\\
 e_N^2(\pi_{opt}) &\sim&    C_3 \, N^{-13/10}         & \ninf,
\end{array}
$$
where $C_1\simeq 0.26$, $C_2\simeq0.20$, and $C_3\simeq0.48$. Observe that utilizing the asymptotically optimal intercomponent grid distribution leads to an increased rate of convergence. \medskip

\begin{figure}[htb]
\begin{center}
\begin{tabular}{cc}
\includegraphics[height=1.8in]{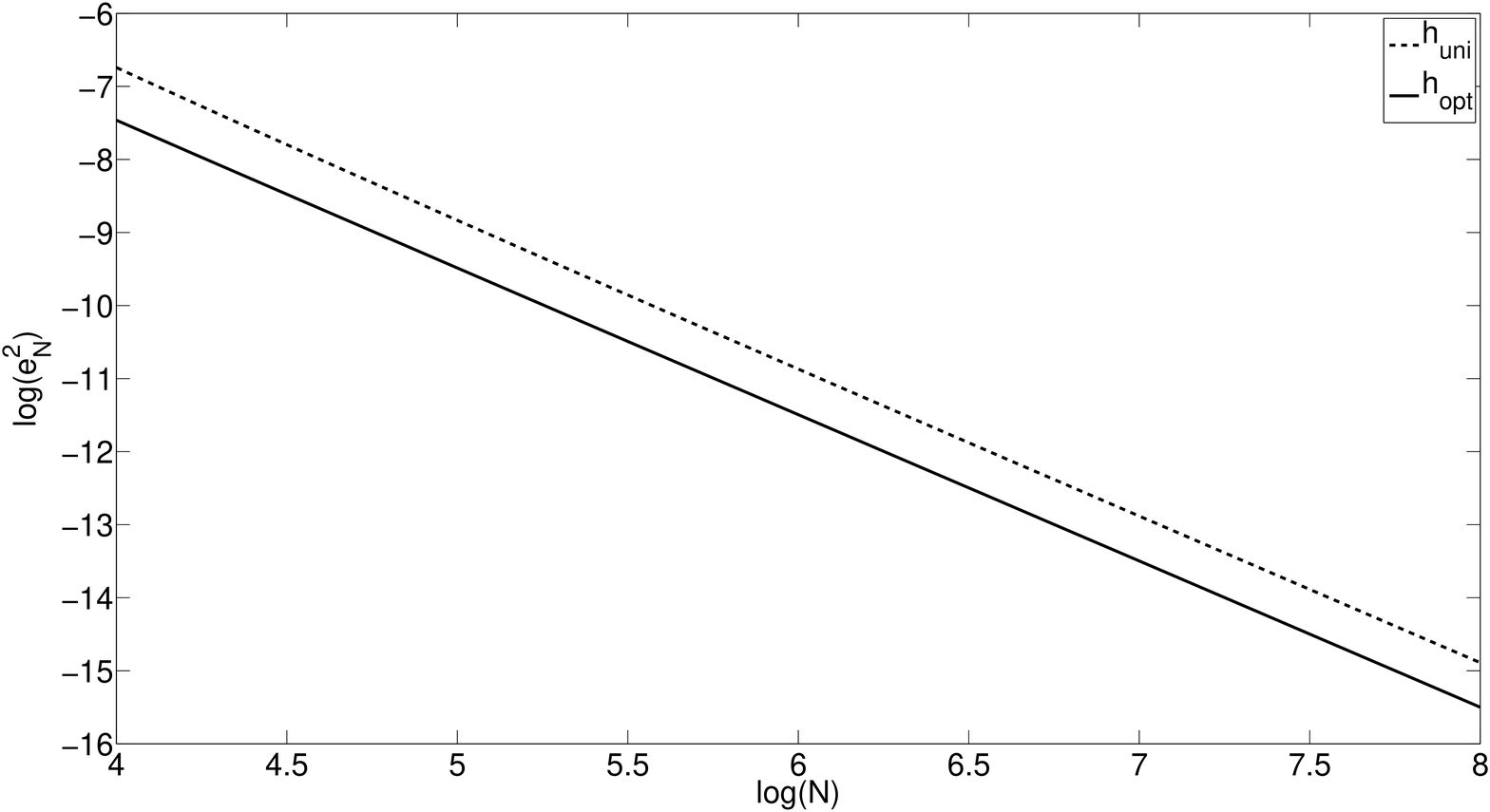} &
\includegraphics[height=1.8in]{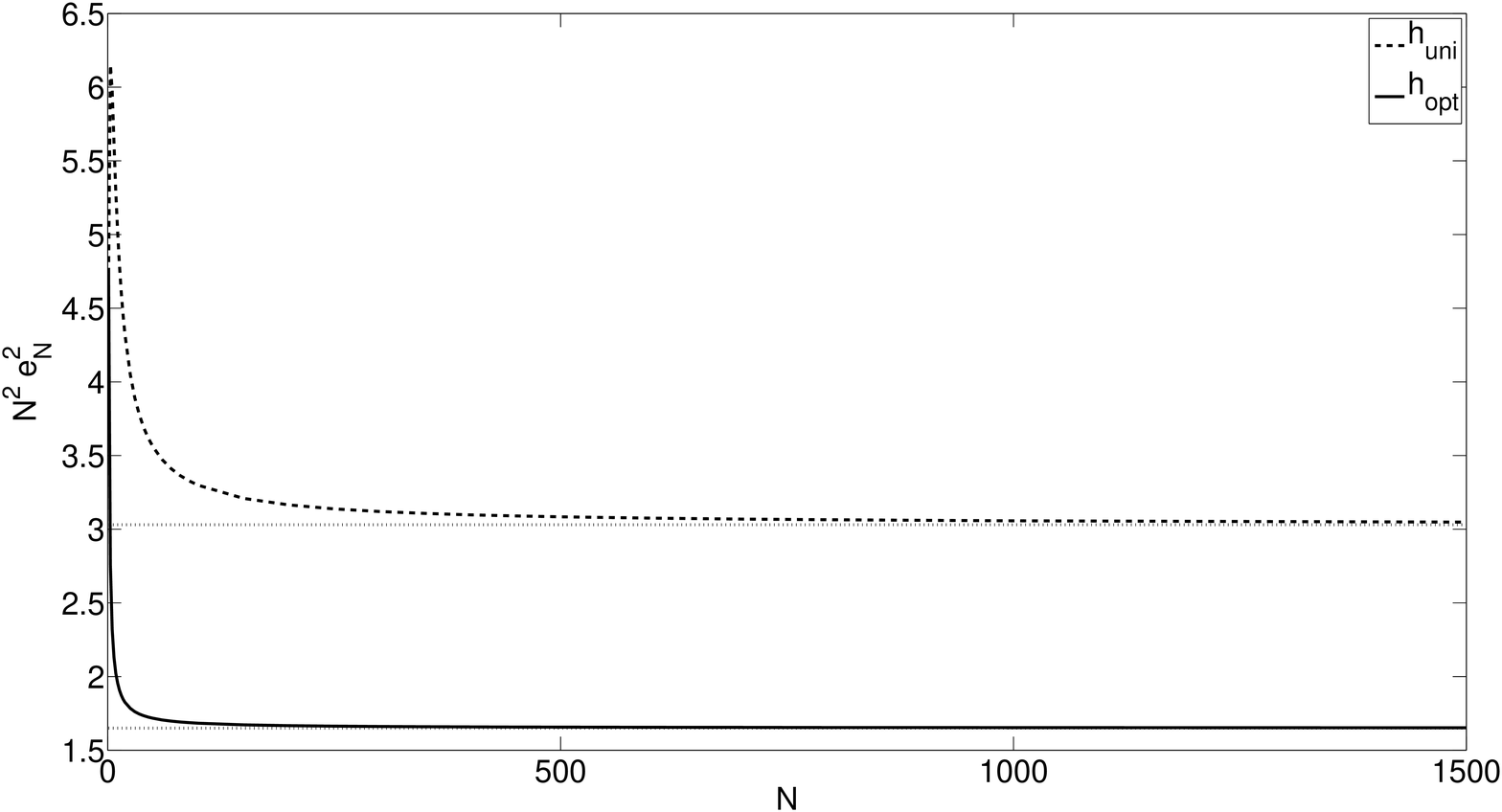}\\
\small{(a)} & \small{(b)}
\end{tabular}
\end{center}
\caption{(a) The (fitted) plots of $e_N^2(h_{uni})$ (dashed line) and $e_N^2(h_{opt})$ (solid line) versus $N$ in a log-log scale.
(b) The convergence of $N^{2}e_N^2(h_{uni})$ (dashed line) and $N^{2}e_N^2(h_{opt})$ to the corresponding asymptotic constants (dotted lines).}
\label{Fig:Ex2}
\end{figure}

\noindent\textbf{Example 4.} Let $Y(t),t\in[0,1]$, be a stochastic process with covariance function $r(t,s)=\exp(-|s-t|)$ and consider process
$$
X(t)=\frac{1}{t+0.1}Y(t),\quad t\in[0,1].
$$
Then $X\in\bcal^\alpha_1([0,1],c(\cdot))$ with $\alpha=1$ and $c(t)=2/(t+0.1)^2$, $t\in[0,1]$.
By Corollary \ref{Cor:MainOneDim}, the squared rate of approximation for any regular density is $N^{-2}$.
We compare the behavior of $e_N^2(h_{uni})=e_N^2(X,h_{uni},\pi_{uni},1)$ and $e_N^2(h_{opt})=e_N^2(X,h_{opt},\pi_{uni},1)$, where $h_{opt}(\cdot)$ given by \eqref{eq:OneDimOpt} is the density minimizing the asymptotic constant.
Figure \ref{Fig:Ex2}(a) shows the (fitted) plots of the mean squared errors $e_N^2(h_{uni})$ (dashed line) and $e_N^2(h_{opt})$ versus $N$ (in a log-log scale). These plots
correspond to the following asymptotic behavior:
$$
\begin{array}{rcl}
 e_N^2(h_{uni}) &\sim&  C_1\, N^{-2},\\
 e_N^2(h_{opt})&\sim&   C_2 \, N^{-2} \ninf
\end{array}
$$
with $C_1\simeq3.03$ and $C_2\simeq1.65$.
Figure \ref{Fig:Ex2}(b) demonstrates the convergence of the scaled  mean squared errors $N^{2}e_N^2(h_{uni})$ and $N^{2}e_N^2(h_{opt})$ to
the corresponding asymptotic constants obtained in Corollary~\ref{Cor:MainOneDim}. Note the
benefit in the asymptotic constant for the optimal density $h_{opt}(\cdot)$.
\medskip

\begin{figure}[htb]
\begin{center}
\includegraphics[height=1.8in]{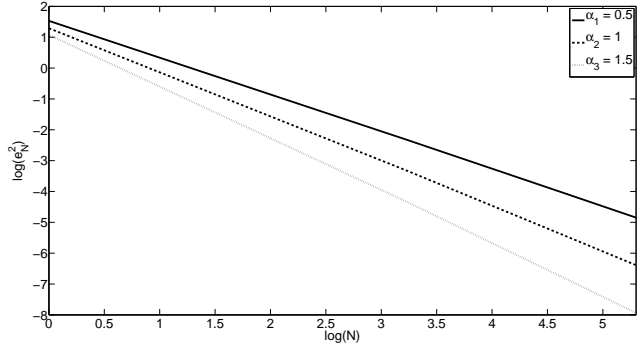}
\end{center}
\caption{The (fitted) plots of $e_N^2(X_{\alpha_j})$, $j=1,2,3$, for $\alpha_1=1/2$ (solid line), $\alpha_2=1$ (dashed line), and $\alpha_3=3/2$ (dotted)  versus $N$ in a log-log scale.}
\label{Fig:Ex3}
\end{figure}

\noindent\textbf{Example 5.} Consider a random field $X_{\alpha}(\mb{t})=\sqrt{10}Y_{\alpha,1/5}$, $\mb{t}\in[0,1]^2$,
where $Y_{\alpha,\beta}$ is defined in the Example 2.
We compare the behavior of the mean squared errors $e_N^2(X_{\alpha_j})=e_N^2(X_{\alpha_j},\pi_{uni},h_{uni})$, $j=1,2,3$, with $\alpha_1=1/2$, $\alpha_2=1$, and $\alpha_3=3/2$.
The local H\"{o}lder function $V(\mb{t})=\norm{\mb{t}}^{-9/5}+2$, $\mb{t}\in[0,1]^2$, satisfies the condition (C).
Consequently by Theorem \ref{Th:RecoveryField},
the sMCQ with cross regular grid sequences attains the convergence rate $N^{-(1+\alpha_j/2)}$, $j=1,2,3$, respectively, despite the point singularity at origin.
Figure \ref{Fig:Ex3} shows the fitted plots of the mean squared errors $e_N^2(X_{\alpha_j})$, $j=1,2,3$ versus $N$ (in a log-log scale).
\medskip

\begin{figure}[htb]
\begin{center}
\begin{tabular}{cc}
\includegraphics[height=1.8in]{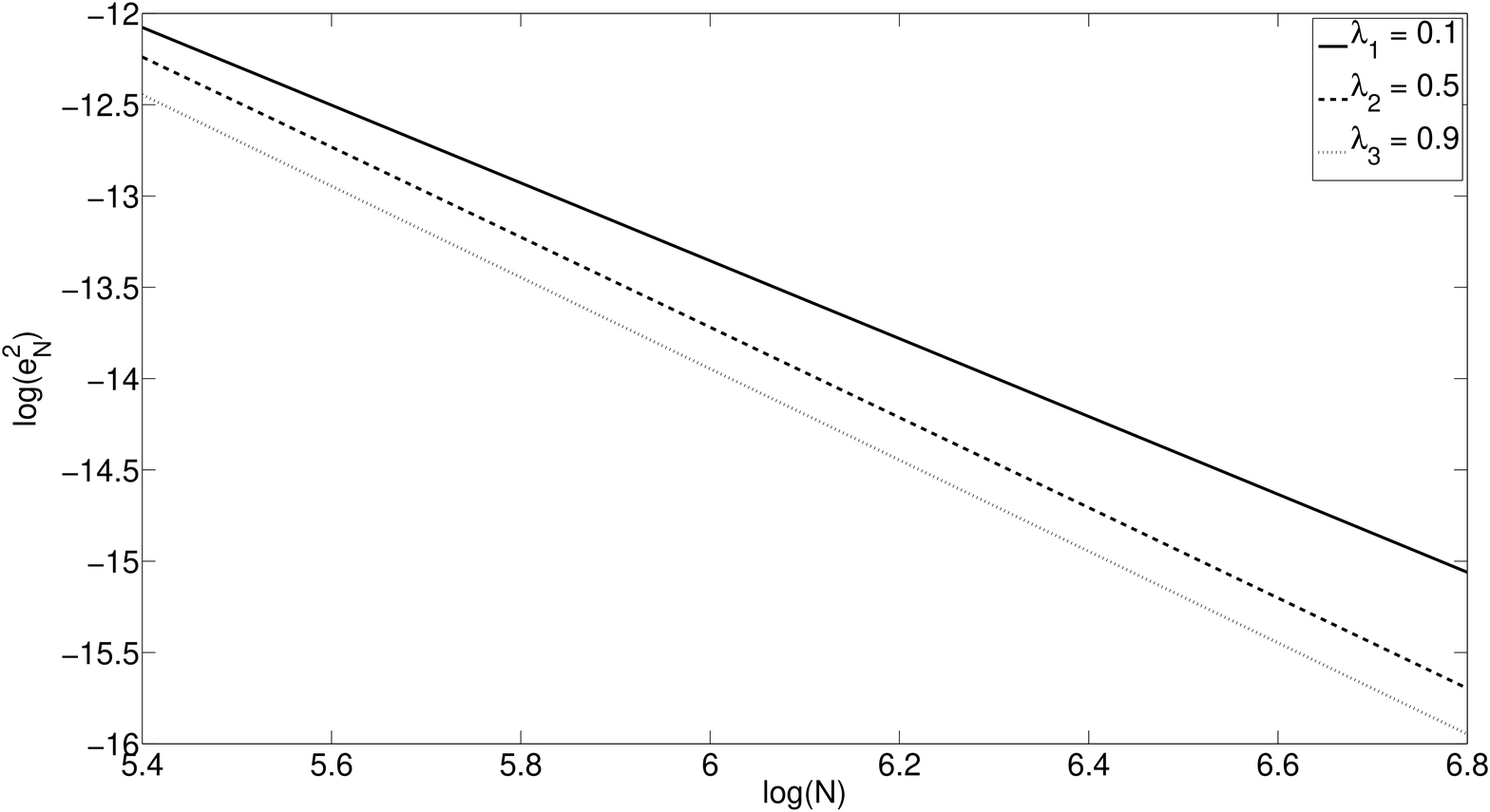} &
\includegraphics[height=1.8in]{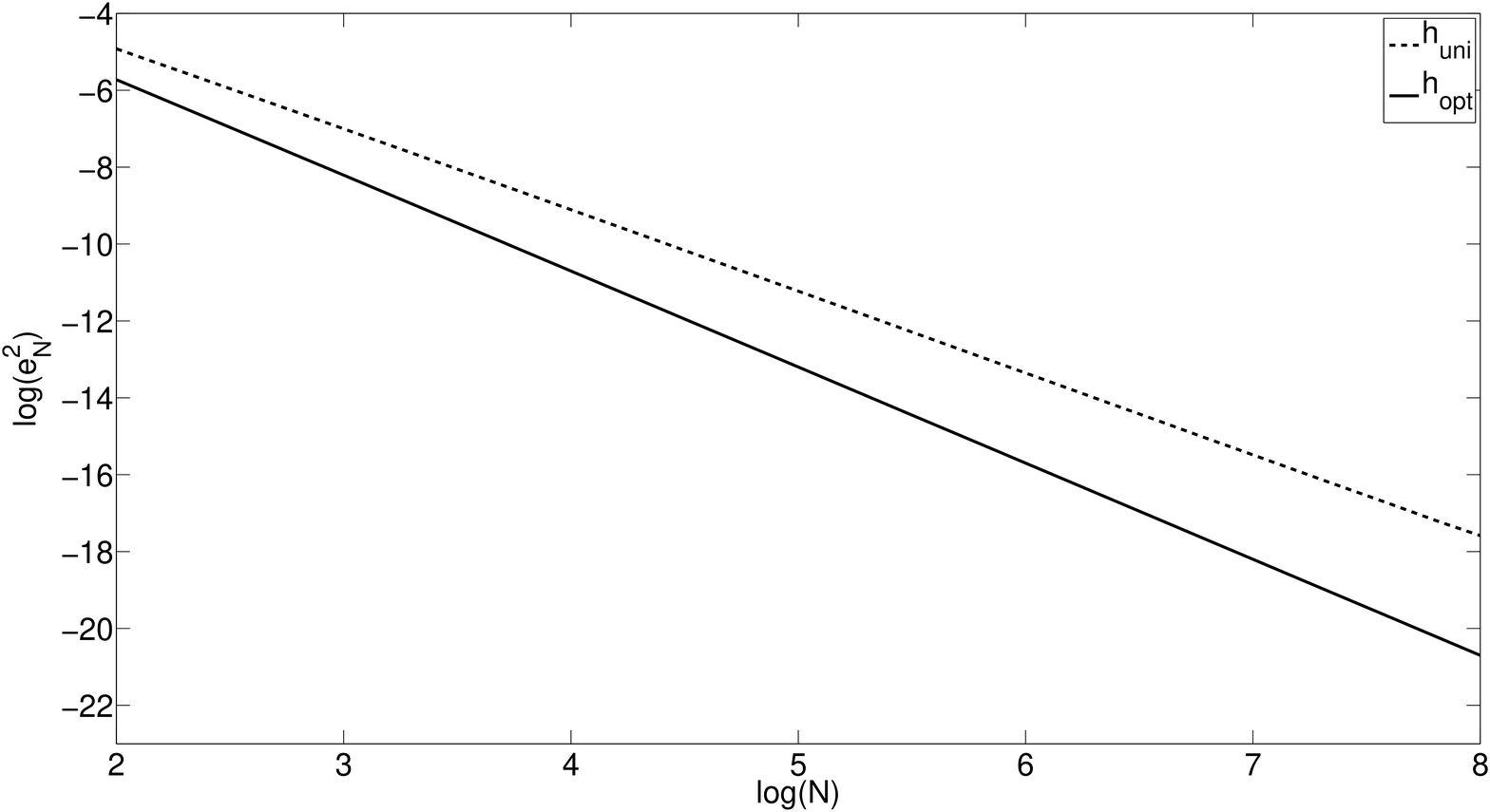}\\
\small{(a)} & \small{(b)}
\end{tabular}
\end{center}
\caption{(a) The (fitted) plots of $e_N^2(X_{\lambda_i},h_{uni})$, $i=1,2,3$ for $\lambda_1=1/10$ (solid line),  $\lambda_2=1/2$ (dashed line) and $\lambda_3=9/10$ (dotted line) versus $N$ in a log-log scale.
(b) The (fitted) plots of $e_N^2(X_{\lambda_1},h_{uni})$ (dashed line) and $e_N^2(X_{\lambda_1},h_{opt})$ (solid line) versus $N$ in a log-log scale.}
\label{Fig:Ex4}
\end{figure}

\noindent\textbf{Example 6.} Let $X_\lambda(t)=5 B_{1,3/2}(t^\lambda)$, $t\in[0,1]$, $0<\lambda<1$, where $B_{m,\beta}$ is a fractional Brownian motion with the covariance function \eqref{eq:FBFdef}.
Then
$$
X_\lambda\in\ccal_1^{3/2\lambda}([0,1],M)\cap\ccal\bcal_1^{3/2}((0,1],c(\cdot),V(\cdot))
$$
with $M=5$ and $c(t)=V(t)=25\lambda^{3/2} t^{3/2(\lambda-1)}$, $t\in[0,1]$. We consider the behavior of the  mean squared errors for $\lambda_1=1/10$, $\lambda_2=1/2$, and $\lambda_3=9/10$.
By Theorem \ref{Th:RecoveryField},
we know that sMCQ with regular grid sequences attains the optimal rate of convergence in two latter cases.
Figure \ref{Fig:Ex4}(a) presents the fitted plots of $e_N^2(X_{\lambda_j},h_{uni})=e_N^2(X_{\lambda_j},h_{uni},\pi_{uni},d)$, $j=1,2,3$.
These plots correspond to the following asymptotic behavior:
$$
\begin{array}{rcl}
 e_N^2(X_{\lambda_1},h_{uni}) &\sim&  C_1\, N^{-2.15},\\
 e_N^2(X_{\lambda_2},h_{uni}) &\sim&  C_2\, N^{-2.5}, \\
 e_N^2(X_{\lambda_3},h_{uni}) &\sim&  C_3\, N^{-2.5}  \ninf
\end{array}
$$
with $C_1\simeq0.64$, $C_2=3.69$, and $C_3\simeq 2.86$.
Consider now the case $\lambda_1=1/10$. By Corollary \ref{Cor:MainOneDim}, the density minimizing the asymptotic constant is given by \eqref{eq:OneDimOpt}.
Moreover, for such defined $h_{opt}(\cdot)$ the condition (C$'$) is satisfied and by Theorem~\ref{Th:Recovery}, the corresponding convergence rate is $N^{-2.5}$.
Figure \ref{Fig:Ex4}(b) shows the (fitted) plots of $e_N^2(X_{\lambda_1},h_{uni})$ and $e_N^2(X_{\lambda_1},h_{opt})$ versus $N$ in a log-log scale. These plots
correspond to the following asymptotic behavior:
$$
\begin{array}{rcl}
 e_N^2(X_{\lambda_1},h_{uni}) &\sim&  C_1\, N^{-2.15},\\
 e_N^2(X_{\lambda_1},h_{opt}) &\sim&  C_4\, N^{-2.5}  \ninf
\end{array}
$$
with $C_4\simeq 0.49$ and an increasing convergence rate for the asymptotically optimal density.

\section{Proofs}
\noindent\textit{Proof of Theorem \ref{Th:Main}.}
Let us recall the definitions:
$$
I(X)=\int_\dcal X(\mb{t}) d\mb{t},
\qquad
I_N(X)=\sum_{\mb{i}\in\mb{I}}  X(\boldsymbol \eta_{\mb{i}})|\dcal_i|,
$$
where $\boldsymbol\eta_{\mb{i}}$ is uniformly distributed in the hyperrectangle $\dcal_{\mb{i}}$, $\mb{i}\in\mb{I}$.
Define the error of numerical integration
$$
\delta_N(X):=I(X)-I_N(X)=\int_\dcal X(\mb{t})d\mb{t}- \sum_{\mb{i}\in\mb{I}}  X(\boldsymbol \eta_{\mb{i}})|\dcal_\mb{i}|=
\sum_{\mb{i}\in\mb{I}}\int_{\dcal_\mb{i}}(X(t)-X(\boldsymbol \eta_{\mb{i}}))d\mb{t},
$$
where  $\E_{\boldsymbol\eta} \delta_N(X)=0$. Denote by
$ e_N^2=e_N^2(X): = \E\delta_N(X)^2 $ the corresponding mean squared error.
By the uniformity and independence of $\boldsymbol \eta_\mb{i}, \mb{i}\in\mb{I}$, we obtain the following expression for the MSE:
\begin{eqnarray}
e_N^2&=& \E\delta_N(X)^2=\E_X\E_{\boldsymbol\eta}\delta_N(X)^2
     =\E_X\V_{\boldsymbol\eta}(\delta_N(X))
     =\E_X\left(\sum_{\mb{i}\in\mb{I}}\V_{\boldsymbol\eta}\left(\int_{\dcal_{\mb{i}}}(X(\mb{t})-X(\boldsymbol \eta_{\mb{i}}))d\mb{t}\right)\right)\nonumber\\
     &=&
     \E_X\left(\sum_{\mb{i}\in\mb{I}}\E_{\boldsymbol\eta}\left(\int_{\dcal_{\mb{i}}}(X(\mb{t})-X(\boldsymbol \eta_{\mb{i}}))d\mb{t}\right)^2\right)
     =\E_X\left(\sum_{\mb{i}\in\mb{I}}\E_{\boldsymbol\eta}\left(\int_{\dcal_{\mb{i}}}\int_{\dcal_{\mb{i}}}(X(\mb{t})-X(\boldsymbol \eta_{\mb{i}}))(X(\mb{s})-X(\boldsymbol \eta_{\mb{i}}))d\mb{t}d\mb{s}\right)\right)\nonumber\\
     &=&\frac{1}{2}\sum_{\mb{i}\in\mb{I}}\E_{\boldsymbol\eta}\left(\int_{\dcal_{\mb{i}}}\int_{\dcal_{\mb{i}}}\left(d_X(\mb{t},\boldsymbol \eta_{\mb{i}})+d_X(\mb{s},\boldsymbol \eta_{\mb{i}})-d_X(\mb{s},\mb{t})\right)d\mb{t}d\mb{s}\right)
     =\frac{1}{2}\sum_{\mb{i}\in\mb{I}}\int_{\dcal_{\mb{i}}}\int_{\dcal_{\mb{i}}}d_X(\mb{t},\mb{v})d\mb{t}d\mb{v},
    \label{eq:ProofMain1}
\end{eqnarray}
where $d_X(\mb{s},\mb{t}):=||X(\mb{t})-X(\mb{s})||^2$ is the incremental variance of the random field $X$. Now the local stationarity condition \eqref{def:locstat} implies that
\begin{equation}\label{eq:Th1proof1}
\begin{aligned}
e_N^2=\frac{1}{2}\left(\sum_{\mb{i}\in\mb{I}}\sum_{j=1}^k c_j(\mb{t_i})
\int_{\dcal_{\mb{i}}}\int_{\dcal_{\mb{i}}} 
\left(\norm{\mb{t}^j-\mb{v}^j}^{\alpha_j}\right)
d\mb{t}d\mb{v}\right)(1+q_{N,\mb{i}}),
\end{aligned}
\end{equation}
where by the positiveness and uniform continuity of local stationarity functions, we have that  $\varepsilon_N=\max\{|q_{N,\mb{i}}|,\mb{i}\in\mb{I}\}=\mathrm{o}(1)$
as $N\to\infty$ \citep[cf.][]{AbramowiczSeleznjev2011}. Recall that the hyperrectangle $\dcal_\mb{i}$ is determined by the vertex
\mbox{$\mathbf{t_{i}}=(t_{1,i_1},\ldots,t_{d,i_d})$}
and the main diagonal $\mathbf{r}_{\mathbf{i}}=\mathbf{t}_{\mb{i}+\mb{1}_d}-\mathbf{t_{i}}$, i.e.,
$$
\dcal_{\mathbf{i}}:=\left\{\mathbf{t}: \mathbf{t}=\mathbf{t}_{\mathbf{i}}+\mathbf{r}_{\mathbf{i}} *\mathbf{s}, \mathbf{s}\in[0,1]^d\right\}.
$$
It follows from the definition and the mean (integral) value theorem that
\begin{equation}\label{eq:rdef}
\mathbf{r}_{\mathbf{i}}=(r_{1,i_1},\ldots,r_{d,i_d})=\left(\frac{1}{h_{1}^*(w_{1,i_1})n_{1}^*},\ldots,\frac{1}{h_{d}^*(w_{d,i_d})n_{d}^*}\right),
      \quad w_{m,i_m}\in [t_{m,i_m},t_{m,i_m+1}],\,m=1,\ldots,d.
\end{equation}
Denote by $\mathbf{w}_{\mathbf{i}}:=(w_{1,i_1},\ldots,w_{d,i_d})$. By the definition of $cRS(h,\pi,\mb{l})$, we get
$$
\mathbf{r}_{\mathbf{i}}^j=\left(\frac{1}{n_j h_j(w_{L_{j-1}+1,i_{L_{j-1}+1}})},\ldots,\frac{1}{n_j h_j(w_{L_j,i_{L_j}})}\right)=\frac{1}{n_j} D_j(\mb{w}_{\mb{i}}^j),\quad j=1,\ldots,k,
$$
where $D_j(\mb{t}^j)=(1/h_j(t_{L_{j-1}+1}),\ldots,1/h_j(t_{L_{j}}))$, $j=1,\ldots,k$. Consequently, changing variables
$\mb{t}^j=\mb{t}_\mb{i}^j+\overline{\mb{t}}^j\ast\mb{r}_{\mb{i}}^j$,
$\mb{v}^j=\mb{t}_\mb{i}^j+\overline{\mb{v}}^j\ast\mb{r}_{\mb{i}}^j$,
$j=1,\ldots,k$, $\mb{i}\in\mb{I}$, gives
$$
\begin{aligned}
e_N^2&=\frac{1}{2}\Bigg(\sum_{\mb{i}\in\mb{I}}|\dcal_\mb{i}|^{2}\sum_{j=1}^k c_j(\mb{t_i})n_j^{-\alpha_j}
\int_{\dcal^j}\int_{\dcal^j}
\norm{D_j(\mb{w}_\mb{i}^j)\ast(\overline{\mb{t}}^j-\overline{\mb{v}}^j)}^{\alpha_j}
d\overline{\mb{t}}^j d\overline{\mb{v}}^j\Bigg)(1+\mathrm{o}(1))\\
\end{aligned}
$$
$\ninf$. Applying the uniform continuity of withincomponent densities, we obtain that
$$
\begin{aligned}
e_N^2&=\frac{1}{2}\Bigg(\sum_{\mb{i}\in\mb{I}}|\dcal_\mb{i}|^{2}\sum_{j=1}^k c_j(\mb{t_i})n_j^{-\alpha_j}
\int_{\dcal^j}\int_{\dcal^j}
\norm{D_j(\mb{t}_\mb{i}^j)\ast(\overline{\mb{t}}^j-\overline{\mb{v}}^j)}^{\alpha_j}
d\overline{\mb{t}}^jd\overline{\mb{v}}^j\Bigg)(1+\mathrm{o}(1))\\
& = \Bigg(\sum_{\mb{i}\in\mb{I}}|\dcal_\mb{i}|^{2}\sum_{j=1}^k c_j(\mb{t_i})n_j^{-\alpha_j}
b_{\alpha_j,l_j}(D_j(\mb{t_i}))\Bigg)(1+\mathrm{o}(1)) \ninf,
\end{aligned}
$$
where $b_{\alpha_j,l_j}(\cdot)$, $j=1,\ldots,k$, are defined by \eqref{def:bconst}.
By equation \eqref{eq:rdef}, we have that
$$|\dcal_{\mb{i}}|= \prod_{m=1}^{d}\frac{1}{n_m^* h^*_m(w_{m,i_m})} =\frac{1}{N} \prod_{m=1}^{d}\frac{1}{h^*_m(w_{m,i_m})}$$
with $w_{m,i_m}\in[t_{m,i_m},t_{m,i_m+1}]$, $\mb{i}\in\mb{I}$, $m=1,\ldots,d$. Furthermore, the uniform continuity of the withincomponent densities implies
$$
\begin{aligned}
e_N^2&=\Bigg(\sum_{\mb{i}\in\mb{I}}|\dcal_\mb{i}| \frac{1}{N}
\prod_{m=1}^{d}\frac{1}{h^*_m(t_{m,i_m})}
 \sum_{j=1}^k c_j(\mb{t_i})n_j^{-\alpha_j}
b_{\alpha_j,l_j}(D(\mb{t_i}^j))\Bigg)(1+\mathrm{o}(1))\\
&=\Bigg(\frac{1}{N} \sum_{j=1}^k  n_j^{-\alpha_j}
\sum_{\mb{i}\in\mb{I}}
 c_j(\mb{t_i}) b_{\alpha_j,l_j}(D(\mb{t_i}^j))
\prod_{m=1}^{d}
{h_m^*(t_{m,i_m})}^{-1}  |\dcal_\mb{i}| \Bigg)(1+\mathrm{o}(1))\ninf.
\end{aligned}
$$
Finally, the Riemann integrability of
$c_j(\mb{t})b_{\alpha_j,l_j}(D(\mb{t}))\prod_{m=1}^{d}{h^*_m(t_m)}^{-1}$
gives
$$
\begin{aligned}
e_N^2&=
\Bigg(\frac{1}{N} \sum_{j=1}^k  n_j^{-\alpha_j}
\int_\dcal
 c_j(\mb{t}) b_{\alpha_j,l_j}(D(\mb{t}^j))
\prod_{m=1}^{d}
{h^*_m(t_m)}^{-1}d\mb{t}\Bigg)(1+\mathrm{o}(1))=\Bigg(\frac{1}{N} \sum_{j=1}^k \frac{v_j}{n_j^{\alpha_j}}\Bigg)(1+\mathrm{o}(1))
\end{aligned}
$$
$\ninf$. This completes the proof.\bigskip

\noindent{\it Proof of Theorem \ref{Th:DimOpt}.} The proof is based on the inequality for the arithmetic and geometric means (cf. Abramowicz and Seleznjev, 2011a), i.e.,
$$
\frac{1}{k}\sum_{j=1}^{k}\frac{v_j}{n_j^{\alpha_j}}\geq\left(\prod_{j=1}^{k}\frac{v_j}{n_j^{\alpha_j}}\right)^{1/k}
$$
with equality if only if
$$
\nu^{-1}=\frac{v_j}{n_j^{\alpha_j}},\quad j=1,\ldots,k.
$$
Hence, the equality is attained for $\tilde n_j=(\nu v_j)^{1/\alpha_j}$, $j=1,\ldots,k$. Let
\begin{equation}\label{eq:ThDimOpt1}
n_j = \left\lceil\tilde n_j\right\rceil\sim\left({\nu}{v_j}\right)^{1/\alpha_j} \mbox{ as } N\to\infty.
\end{equation}
This implies that for the asymptotically optimal intercomponent knot distribution
$$
N\sim \nu^{1/\rho}\prod_{j=1}^{k}{v_j^{l_j/\alpha_j}},
$$
and therefore,
$$
\nu\sim {N^{\rho}}{\kappa^{-\rho}} \mbox{ as } N\to\infty.
$$
By equation \eqref{eq:ThDimOpt1}, the asymptotically optimal intercomponent knot distribution is
$$
n_j\sim\frac{N^{\rho/\alpha_j} v_j^{1/\alpha_j} }{\kappa^{\rho/\alpha_j}}\mbox{ as } N\to\infty,\quad j=1,\ldots,k.
$$
Moreover, with such knot distribution, the equality in \eqref{eq:Th2} is attained asymptotically. This completes the proof.\bigskip

\noindent{\it Proof of Proposition \ref{Prop:OneDim}.} The proof follows directly from the proof of Theorem~1. The expression for
the optimal withincomponent density follows from Seleznjev(2000).
\bigskip

\noindent{\it Proof of Proposition \ref{Prop:Holder}.} The first steps of the proof repeat those of Theorem \ref{Th:Main}. Applying the H\"{o}lder condition \eqref{def:hcont} to equation \eqref{eq:ProofMain1}
yields
$$
\begin{aligned}
e_N^2
&\leq \frac{1}{2} C \sum_{\mb{i}\in\mb{I}} \sum_{j=1}^k \int_{\dcal_{\mb{i}}}\int_{\dcal_{\mb{i}}} ||\mb{t}^j-\mb{v}^j||^{\alpha_j}d\mb{t}d\mb{v}
\leq \frac{1}{2}  C \sum_{\mb{i}\in\mb{I}} \sum_{j=1}^k l_j^{\alpha_j/2}\sum_{m=L_{j-1}+1}^{L_j} \int_{\dcal_{\mb{i}}}\int_{\dcal_{\mb{i}}} |t_m-v_m|^{\alpha_j}d\mb{t}d\mb{v},
\end{aligned}
$$
where the last inequality follows from the fact that any nonnegative numbers $a_1,\ldots,a_k$ and any $\alpha\in \mathbb{R_+}$, the inequality
\begin{equation}\label{PowerIneq}
\left(\sum_{i=1}^{k} a_i\right)^\alpha \leq k^\alpha\sum_{i=1}^{k}a_i^\alpha
\end{equation}
holds. Consequently, changing variables $\bar t=(t_m-t_{m,i_m})/r_{m,i_m}$, $\bar v=(v_m-t_{m,i_m})/r_{m,i_m}$, $m=1,\ldots,d$, $\mb{i}\in\mb{I}$, gives
$$
\begin{aligned}
e_N^2& \leq \frac{1}{2} C \sum_{\mb{i}\in\mb{I}} \sum_{j=1}^k l_j^{\alpha_j/2} |\dcal_\mb{i}|^2\sum_{m=L_{j-1}+1}^{L_j} r_{m,i_m}^{\alpha_j} \int_{0}^1\int_{0}^1 |\bar t-\bar v|^{\alpha_j}d{\bar t}d{\bar v}
= C \sum_{\mb{i}\in\mb{I}} \sum_{j=1}^k l_j^{\alpha_j/2} a_{\alpha_j} |\dcal_\mb{i}|^2\sum_{m=L_{j-1}+1}^{L_j} r_{m,i_m}^{\alpha_j},
\end{aligned}
$$
where $a_\alpha=1/2\int_0^1\int_0^1 |t-s|^\alpha dt ds=1/((1+\alpha)(2+\alpha))$.
By the continuity of the withincomponent densities and the mean value theorem, we have that
\begin{equation}\label{eq:rbound}
r_{m,i_m}\leq \frac{C_m^*}{n_m^*},\qquad \mb{i}\in\mb{I},\, m=1,\ldots,d
\end{equation}
with $C_m^*=1/\min_{s\in[0,1]}h^*_m(s)$.
Moreover, the definition of $cRS(h,\pi,\mb{l})$ implies that
$$
\begin{aligned}
e_N^2& \leq  C \left(\prod_{m=1}^{d}\frac{C_m^*}{n_m^*}\right) \sum_{j=1}^k a_{\alpha_j} l_j^{1+\alpha_j/2}  \left(\frac{C^*_{L_j}}{n_j}\right)^{\alpha_j} \sum_{\mb{i}\in\mb{I}}|\dcal_\mb{i}|
=\frac{C}{N} \sum_{j=1}^{k}\frac{d_j}{{n_j}^{\alpha_j}}
\end{aligned}
$$
with $d_j=a_{\alpha_j} l_j^{1+\alpha_j/2}(C_{L_j}^*)^{\alpha_j}\prod_{m=1}^{d}C^*_m$,  $j=1,\ldots,k$. The formula for the asymptotically optimal intercomponent grid distribution follows
from the proof of Theorem \ref{Th:DimOpt}. This completes the proof.\bigskip\\

\noindent{\it Proof of Theorem \ref{Th:RecoveryField}.} The first steps of the proof repeat those of Theorem 1.
Consider equation \eqref{eq:ProofMain1}. The MSE can be decomposed as follows:
$$
\begin{aligned}
e_N^2=\frac{1}{2}\sum_{\mb{i}\in\mb{I}}\int_{\dcal_{\mb{i}}}\int_{\dcal_{\mb{i}}}d_X(\mb{t},\mb{v})d\mb{t}d\mb{v}
=\sum_{\mb{i}\in\mb{I}}e_{\mb{i},N}^2
\end{aligned}
$$
with
$$
e_{\mb{i},N}^2= \frac{1}{2} \int_{\dcal_{\mb{i}}}\int_{\dcal_{\mb{i}}}d_X(\mb{t},\mb{v})d\mb{t}d\mb{v}.
$$
For a fixed $\delta>0$, we denote $\Delta:=[0,\delta]^d$, and $\mb{I}_\Delta:=\{\mb{i}:\dcal_\mb{i} \cap \Delta \neq \emptyset\}$. Consequently,
$$
e_N^2=\sum_{\mb{i}\in\mb{I}}e_{\mb{i},N}^2=e_{\mb{0}_d,N}^2 + \sum_{\mb{i}\in{\mb{I}_{\Delta}\backslash\{\mb{0}_d\}}} e_{i,N}^2+ \sum_{\mb{i}\in{\mb{I}\backslash\mb{I}_{\Delta}} } e_{i,N}^2=S_1+S_2+S_3,
$$
where $S_1=S_1(N):=e_{\mb{0}_d,N}^2$,
$S_3=S_3(N)$ includes all terms $e_{\mb{i},N}$ such that $\dcal_\mb{i}\subset \dcal\backslash\Delta$,
 and
$S_2=S_2(N):=e_N^2-S_1-S_3$.
For $S_1$,  the H\"{o}lder condition, \eqref{PowerIneq}, and \eqref{eq:rbound} imply that
\begin{eqnarray}
e_{\mb{0}_d,N}^2
&\leq& C\int_{\dcal_{\mb{0}_d}}\int_{\dcal_{\mb{0}_d}}\norm{\mb{t}-\mb{v}}^{\beta}d\mb{t}d\mb{v}
\leq C d^{\beta/2}\sum_{m=1}^{d} \int_{\dcal_{\mb{0}_d}}\int_{\dcal_{\mb{0}_d}}|t_m-v_m|^{\beta}d\mb{t}d\mb{v}\nonumber\\
&\leq&  C|\dcal_{\mb{0}_d}|^2  d^{\beta/2}a_{\beta}\sum_{m=1}^{d} r_{m,0}^{\beta}
\leq C_1 d^{1+\beta/2}a_{\beta}  {N^{-(2+\beta/d)}}  \label{eq:S_1}
\end{eqnarray}
for a positive constant $C_1$. Hence $
e_{\mb{0}_d,N}^2=o(N^{-(1+\alpha/d)})
$ for any $\beta\in(0,2)$, $\alpha\in(0,2)$, if $d\geq2$, and for $\beta>\alpha-1$, if $d=1$.
For $S_2$ by the local H\"{o}lder condition \eqref{def:locHol}, we obtain the following upper bound
$$
\begin{aligned}
S_2&
=\sum_{\mb{i}\in \mb{I}_\Delta\backslash\{\mb{0}_d\}}e_{\mb{i},N}^2
\leq \sum_{\mb{i}\in \mb{I}_\Delta\backslash\{\mb{0}_d\}} \int_{\dcal_\mb{i}}\int_{\dcal_\mb{i}}d_X(\mb{t},\mb{s})d\mb{t}d\mb{s}
\leq\sum_{\mb{i}\in \mb{I}_\Delta\backslash\{\mb{0}_d\}}  V(\mb{v}_{\mb{i}}) \int_{\dcal_{\mb{i}}}\int_{\dcal_{\mb{i}}}\norm{\mb{t}-\mb{s}}^{\alpha}d\mb{t}d\mb{s}
\end{aligned}
$$
for $\mb{v}_{\mb{i}}\in\dcal_\mb{i}$, $\mb{i}\in \mb{I}_\Delta\backslash\{\mb{0}_d\}$.
The continuity of withincomponent grid generating densities together and the definition of function $\mb{G}(\cdot)$ and condition (C) give
$$
S_2 \leq C_1 {N^{-(1+\alpha/d)}}
\sum_{\mb{i}\in\mb{I}_\Delta\backslash\{\mb{0}_d\}}
 V(\mb{G}(\mb{w}_{\mb{i}}))  |\dcal_\mb{i}|
\leq C_1 {N^{-(1+\alpha/d)}} \sum_{\mb{i}\in\mb{I}_\Delta\backslash\{\mb{0}_d\}}
 R(\mb{w}_{\mb{i}})  |\dcal_\mb{i}|,
$$
where $C_1$ is a positive constant and $\mb{w}_\mb{i}\in[i_1/n_1^*,(i_1+1)/n_1^*]\times\ldots\times[i_d/n_d^*,(i_d+1)/n_d^*]=:\dcal^*_\mb{i}$. The shifting property \eqref{eq:shift} implies that for a positive constant
$C_2$,
\begin{equation}\label{eq:ThRecEq1}
S_2 \leq C_2 {N^{-(1+\alpha/d)}} \sum_{\mb{i}\in\mb{I}_\Delta\backslash\{\mb{0}_d\}}  R(\mb{s}_{\mb{i}})  |\dcal_\mb{i}|,
\end{equation}
where $\mb{s}_\mb{i}=\mb{H}(\mb{u_i})\in\dcal^*_\mb{i}$ is such that
$$R(\mb{s_i})=R(\mb{H}(\mb{u_i}))=\min_{\mb{v_{i}}\in\dcal_\mb{i}}R(\mb{H}(\mb{v}_\mb{i})),\qquad \mb{i}\in \mb{I}_\Delta\backslash\{\mb{0}_d\}.$$
Consequently by \eqref{eq:ThRecEq1} and condition (C), we have
$$
S_2 \leq C_2 {N^{-(1+\alpha/d)}}
\int_{\Delta\backslash\dcal_{\mb{0}_d}} R(\mb{H}(\mb{t}))d\mb{t}.
$$
Thus for any $\varepsilon>0$ and sufficiently small $\delta$ by condition (C), we obtain that
\begin{equation}\label{eq:S2}
N^{1+\alpha/d}S_2\leq C \int_{\Delta\backslash\dcal_{\mb{0}_d}} R(\mb{H}(\mb{t}))d\mb{t} <\varepsilon.
\end{equation}
For $S_3$, similarly to Theorem \ref{Th:Main}, we get that
\begin{equation}\label{eq:S3_1}
 N^{1+\alpha/d}S_3=\left(\int_{\dcal\backslash\Delta}c(\mb{t})b_{\alpha,d}(D_1(\mb{t}))\prod_{m=1}^d h(t_m)^{-1}d\mb{t}\right)(1+\mathrm{o}(1)):=v_{1,\delta}(1+\mathrm{o}(1)) \ninf,
\end{equation}
where $D_1(\mb{t})=(1/h(t_1),\ldots,1/h(t_d))$.
From the regularity of the withincomponent density and condition (C) it follows that for a positive constant $C_1$,
$$\int_\dcal c(\mb{t})b_{\alpha,d}(D_1(\mb{t}))\prod_{m=1}^d h(t_m)^{-1}d\mb{t}\leq C_1 \int_\dcal c(\mb{t}) d\mb{t} \leq C_1 \int_\dcal V(\mb{t}) d\mb{t} \leq C_1 \int_\dcal R(H(\mb{t}))d\mb{t} < \infty$$
and therefore the monotone convergence gives
\begin{equation}\label{eq:S3_2}
v_{1,\delta}\uparrow v_1=\int_{\dcal}c(\mb{t})b_{\alpha,d}(D_1(\mb{t}))\prod_{m=1}^d h(t_m)^{-1}d\mb{t}\quad\mbox{as }\delta\to 0.
\end{equation}
So, for any $\varepsilon>0$, first we select $\delta$ sufficiently small and apply \eqref{eq:S2} and \eqref{eq:S3_2}. Then for the selected $\delta$ and sufficiently large~$N$,
\eqref{eq:S_1} and \eqref{eq:S3_1} imply the assertion.
 This completes the proof. \medskip

\noindent{\it Proof of Theorem \ref{Th:Recovery}.}
The first steps of the proof repeat those of Theorem 1. Consider equation \eqref{eq:ProofMain1} and decompose the MSE as in the proof of Theorem \ref{Th:RecoveryField}:
$$
e_N^2=\sum_{i=0}^{N-1} \frac{1}{2} \int_{t_{i}}^{t_{i+1}}\int_{t_{i}}^{t_{i+1}}d_X({t},{v})d{t}d{v}
=\sum_{i=0}^{N-1}e_{i,N}^2
$$
with
$$
e_{i,N}^2=\frac{1}{2}\int_{t_{i}}^{t_{i+1}}\int_{t_{i}}^{t_{i+1}}d_X({t},{v})d{t}d{v},\quad i=1,\ldots,N.
$$
Moreover, let for fixed $\delta>0$
$$
e_N^2=\sum_{i=0}^{N-1}e_{i,N}^2=e_{0,N}^2 + \sum_{i=1}^{J_\dt} e_{i,N}^2+ \sum_{j=J_\dt+1}^{N-1} e_{i,N}^2=S_1+S_2+S_3,
$$
where $S_1=S_1(N):=e_{0,N}^2$,  $S_3=S_3(N)$ includes all terms $e_{i,N}$ such that $[t_{i},t_{i+1}]\subset [\dt, 1]$,
say, $i\ge J_\dt+1$, and $S_2=S_2(N):=e_N^2-S_1-S_3$.
For $S_1$, the H\"{o}lder condition and the definition of function $G(\cdot)$ implies that
$$
S_1=\frac{1}{2}\int_{0}^{t_1} \int_{0}^{t_1}d_X(t,v)dtdv\leq \frac{1}{2}M \int_{0}^{t_1} \int_{0}^{t_1}|t-v|^{\beta}dtdv = M t_1^{2+\beta} a_{\beta}
   \leq C G\left(\frac{1}{N}\right)^{2+\beta},
$$
for a positive constant $C$. By  condition (C$'$), we obtain that
\begin{equation}\label{eq:ProcRecEq1}
N^{1+\alpha}S_1\leq C N^{1+\alpha} G\left(\frac{1}{N}\right)^{2+\beta}= C N^{1+\alpha}\mathrm{o}(N^{-(1+\alpha)})=\mathrm{o}(1).
\end{equation}
We proceed to calculating the upper bound for $S_2$. By the local H\"{o}lder continuity \eqref{def:locHol} and the mean value theorem,  we obtain that
$$
\begin{aligned}
S_2=\sum_{i=1}^{J_\dt} e_{i,N}^2&\leq\frac{1}{2} \sum_{i=1}^{J_\delta}\int_{t_{i}}^{t_{i+1}} \int_{t_{i}}^{t_{i+1}}d_X(t,v)dtdv
   \leq\frac{1}{2} \sum_{i=1}^{J_\dt} V(G(w_i)) \int_{t_{i}}^{t_{i+1}} \int_{t_{i}}^{t_{i+1}}|t-v|^{\alpha}dtdv = a_{\alpha} \sum_{i=1}^{J_\dt} V(G(w_i)) r_i^{2+\alpha_j} \\
&\leq C N^{-(1+\alpha)} \sum_{i=1}^{J_\dt} r_i V(G(w_i)) g(v_i)^{1+\alpha},
\end{aligned}
$$
where $w_i,v_i\in[i/N,(i+1)/N]$ and $C$ is a positive constant. Now applying the shifting property \eqref{eq:shift} and condition (C$'$), we get
$$
N^{1+\alpha}S_2 \leq C \sum_{i=1}^{J_\dt} r_i V(G(w_i)) g(v_i)^{1+\alpha}\leq C_1 \sum_{i=1}^{J_\dt} r_i R(s_i) r(s_i)^{1+\alpha}
\leq C_1 \int_{t_1}^{\delta}R(H(t))r(H(t))^{1+\alpha}dt,
$$
where for $s_i=H(u_i)\in[i/N,(i+1)/N]$,
$$
R(s_i)r(s_i)^{1+\alpha}=R(H(u_i))r(H(u_i))^{1+\alpha}=\min_{t\in[t_i,t_{i+1}]}R(H(t))r(H(t))^{1+\alpha}.
$$
Thus for any $\epsilon>0$ and sufficiently small $\delta$ by condition (C$'$), we have
\begin{equation}\label{eq:ProcRecEq2}
N^{1+\alpha}S_2\leq C_1 \int_{t_1}^{\delta}R(H(t))r(H(t))^{1+\alpha}dt<\epsilon.
\end{equation}
For $S_3$, we obtain that
\begin{equation}\label{eq:ProcRecEq3}
N^{1+\alpha}S_3 = a_{\alpha}\int_\delta^1 c(t)h^{-(1+\alpha)}dt (1+o(1))=:q_\delta (1+o(1))\ninf.
\end{equation}
It follows by the equation \eqref{def:qdens} and condition (C$'$) that
$$\int_0^1 c(t)h(t)^{-(1+\alpha)} dt = \int_0^1 c(t)g(H(t))^{1+\alpha} dt \leq \int_0^1 V(t)g(H(t))^{1+\alpha} dt \leq \int_0^1 R(H(t))g(H(t))^{1+\alpha} dt <  \infty$$
and the monotone convergence gives
\begin{equation}\label{eq:ProcRecEq4}
q_\delta \uparrow q := a_\alpha \int_0^1 c(t)h(t)^{-(1+\alpha)}dt\quad\mbox{ as } \delta\to 0.
\end{equation}
So, for any $\varepsilon>0$, first we select $\delta$ sufficiently small and apply \eqref{eq:ProcRecEq2} and \eqref{eq:ProcRecEq4}. Then for the selected $\delta$ and sufficiently large~$N$,
\eqref{eq:ProcRecEq1} and \eqref{eq:ProcRecEq3} imply the assertion.
This completes the proof.
\medskip\hfill

\noindent{\bf Acknowledgments}
\smallskip\hfill

\noindent
The second author is partly supported by the Swedish Research Council grant 2009-4489 and the project "Digital Zoo" funded by the European Regional Development Fund.

\bibliographystyle{model2-names}

\bibliography{IntegrationRandomFields}
\end{document}